\documentstyle[12pt]{article}
\newcommand{\const}{\mathop{\rm const}\limits}

\newcommand{\mod}{\mathop{\rm mod}\limits}

\newcommand{\dom}{\mathop{\rm dom}\limits}

\newcommand{\diam}{\mathop{\rm diam}\limits}

\begin{document}

\begin{center}

{\bf DISCOVERING OF BOUNDEDNESS AND CONTINUITY }\\

\vspace{3mm}

{\bf OF RANDOM FIELDS BY MEANS OF}\\

\vspace{3mm}

{\bf PARTITION ENTROPIC SCHEME.  } \\

\vspace{4mm}

                {\bf Eugene Ostrovsky, Leonid Sirota.}  \\

\vspace{3mm}

Department of  Mathematics, Bar-Ilan University, \\

\vspace{3mm}

 Ramat Gan, 52900, Israel, \\

\vspace{3mm}

 e-mails:\ eugostrovsky@list.ru; \\

 \vspace{3mm}

  sirota3@bezeqint.net. \\

\vspace{4mm}
                    {\sc Abstract.}\\

 \end{center}

 \vspace{3mm}

 \ We construct a new sufficient conditions for boundedness or continuity of arbitrary random fields relying on the
so-called  partition  scheme, alike in the classical majorizing measure method. \par
 \ We deduce also the used in the practice (statistics, method Monte-Carlo etc.) exact exponential estimates
for tail of distribution of maximum for random field satisfying  formulated in this report conditions. \par

 \vspace{3mm}


\vspace{3mm}


\vspace{3mm}

\section{ Introduction. Notations. Statement of problem.} \par

\vspace{3mm}

  \ Let $ (T,d) $ be arbitrary compact metric space  with a set $  T  $ equipped with a distance (or, more generally,
semi-distance)  function $ d = d(t,s), \ t,s \in T. $
Let also $ \xi = \xi(t), \ t \in T $ be numerical valued separable stochastic continuous random process or separable random field. Denote
by the $  {\bf P_c}(\xi) $ a probability that almost every path of $  \xi = \xi(t) $ is continuous:

$$
{\bf P_c}(\xi) \stackrel{def}{=} {\bf P} (\xi(\cdot) \in C(T,d) ), \eqno(1.1_c)
$$
as well as denote by the $  {\bf P_b}(\xi) $ a probability that almost every path of $  \xi = \xi(t) $ is bounded:

$$
{\bf P_b}(\xi) = {\bf P} ( \sup_{t \in T} |\xi(t)| < \infty).  \eqno(1.1_b)
$$

 \  It may be interest in many applications: method Monte-Carlo, statistics etc. the concrete estimate of tail for
maximum  distribution

$$
{\bf P}_T(u) = {\bf P}_{T,\xi}(u) = {\bf P} ( \sup_{t \in T} |\xi(t)| > u), \ u \ge C > 0,
$$
as well as the tail estimate for module of continuity of the r.f. $  \xi(\cdot), $
see  \cite{Frolov1}, \cite{Grigorjeva1}, \cite{Ostrovsky2}, \cite{Ostrovsky3} etc.\par

\vspace{3mm}

 \ {\bf We will talk about the sufficient conditions imposed on the finite-dimensional distributions of the r.f. $  \xi (t) $  under which}
$ {\bf P_c}(\xi) = 1 $ or $ {\bf P_b}(\xi) = 1. $  \par

\vspace{3mm}

 \ This problem may be named as classical. There are at last two approaches for his solving: the so-called entropy approach, where
the metric entropy is calculated  relative the distance, more precisely, semi-distance
 function generated by the r.f. $  \xi(t)  $ itself, see e.g.
\cite{Kolmogorov1},  \cite{Kozachenko1}, \cite{Slutzky1}, \cite{Dudley1},  \cite{Fernique1}, \cite{Fernique2}, \cite{Fernique3},
\cite{Ostrovsky3}, \cite{Ostrovsky7}, \cite{Pizier1}; the second one is named as majorizing measure approach, see \cite{Fernique1},
\cite{Talagrand2} - \cite{Talagrand5}, \cite{Ostrovsky2} etc. \par

\vspace{4mm}

 {\sc  Several notations, definitions and facts.}\par

\vspace{3mm}

 \ {\bf A.} A triplet $ (\Omega, \cal{B}, {\bf P} ),  $  where $ \Omega = \{\omega\} $ or $ \Omega = \{x\}  $ is arbitrary set,
$ \cal{B} $ is non-trivial sigma-algebra of subsets $ \Omega  $ and  $ {\bf P}  $ is non-zero non-negative completely
additive normed:  $ \ {\bf P}(\Omega) = 1 \ $  measure defined on the $  \cal{B} $ is called  as ordinary a probabilistic space.\par

 \ We denote as usually for the random variable $ \xi $ (r.v.) (i.e. measurable function $ \Omega \to R )  $ its
classical Lebesgue-Riesz norm by

 $$
 |\xi|_p = [ {\bf E} |\xi|^p ]^{1/p} = \left[ \int_{ \Omega  } |\xi(\omega)|^p \ {\bf P}(d \omega)  \right]^{1/p}, \ p \ge 1; \eqno(1.2)
 $$

$$
L_p = \{ \xi, \ |\xi|_p < \infty  \}.
$$

\vspace{3mm}

{\bf B.} The so-called Grand Lebesgue Space  $ G \psi $ equipped with a norm $  || \ \cdot \ ||G\psi   $ consists by definition
on all the numerical valued random variables  defined on the our probability space and having a following finite norm

$$
G \psi = \{\xi, \ ||\xi||G \psi < \infty \}, \hspace{5mm}  ||\xi||G \psi \stackrel{def}{=} \sup_{ p \ge 1 } \left[ \frac{|\xi|_p}{\psi(p)} \right].
\eqno(1.3)
$$
 \ Here $ \psi = \psi(p) $ is some continuous strictly positive function such that there exists $ \lim_{p \to \infty} \psi(p) = \infty. $\par
 \ These  spaces are Banach rearrangement invariant complete spaces. The detail investigation of this spaces (and more general ones) see in \cite{Liflyand1},  \cite{Ostrovsky8}. See also \cite{Fiorenza1}, \cite{Fiorenza2}, \cite{Iwaniec1}, \cite{Iwaniec2}, \cite{Kozachenko1} etc. \par

\vspace{4mm}

   {\it An important for us fact about considered here spaces is proved in \cite{Kozachenko1}, \cite{Ostrovsky3}, \cite{Ostrovsky6}:
they coincide with some exponential  Orlicz's spaces}
  $ Or(\Phi_{\psi}). $ For instance, if $ \psi(p) = \psi_{1/2}(p): =\sqrt{p}, \ p \ge 1, $ then the space
  $ G\psi_{1/2} $ consists on all the subgaussian (non-centered, in general case) r.v. $ Or(\Phi_{\psi_{1/2}})  $
   for which by definition $ \Phi_{\psi_{1/2}}(u) = \exp(u^2/2) - 1.  $  \par
 \  The Gaussian distributed r.v. $ \eta $ belongs to this space.  Another example: let $ \Omega = (0,1) $ with usually Lebesgue measure and

 $$
 f_{1/2}(\omega) = \sqrt{|\log \omega}|, \ \omega > 0; \ f_{1/2}(0)= 0.
 $$

  \ It is easy to calculate using Stirling's formula for the Gamma function:

 $$
 | f_{1/2}|_p \asymp \sqrt{p}, \ p \in (1,\infty).
 $$

  \ The correspondent tail behavior:

  $$
  {\bf P} ( f_{1/2} > u ) = \exp (-u^2).
  $$

 \ More detail, let the function $ \psi(\cdot) \in G\Psi = G\Psi_{\infty} $  be such that the new generated by $ \psi $ function

 $$
 \nu(p) = \nu_{\psi}(p) := p \ln \psi(p), \ p \in [1,\infty)
 $$
is convex. The r.v. $ \eta $ belongs to the space $  G\psi $ if and only if it belongs to the Orlicz's space $ L(N_{\psi}) $ with the
correspondent {\it exponential} continuous Young-Orlicz function

$$
N_{\psi}(u) := \exp \left( - \nu^*_{\psi}(\ln |u|)  \right), \ |u| \ge e,
$$

$$
N_{\psi}(u) := C \ u^2, \ |u| < e, \hspace{4mm} C \ e^2 = \exp \left( - \nu^*_{\psi}(1)  \right),
$$
and herewith of course both the Banach space norms: $  ||\cdot|| L(N_{\psi})  $  and $ || \cdot||G\psi  $ are equivalent. \par

 \ One can also complete characterize (under formulated here conditions)  the belonging of the non-zero r.v. $ \xi $ to the
space $ G\psi $  by means of its tail behavior:

$$
\xi \in G\psi \Leftrightarrow \exists K = \const \in (0, \infty), \ \max( {\bf P}(\xi > u), {\bf P}(\xi < - u)  ) \le
$$

$$
\exp \left( - \nu^*_{\psi}(\ln |u|/K)  \right), \ u \ge K e,
$$
see \cite{Kozachenko1}, \cite{Ostrovsky3}, p. 33-35.\par

 \vspace{4mm}

 \  The case when in (1.3) the supremum is calculated over {\it finite } interval is investigated in
 \cite{Liflyand1}, \cite{Ostrovsky8}, \cite{Ostrovsky9}:

$$
G_b \psi = \{\xi, \ ||\xi||G_b \psi < \infty \}, \  ||\xi||G_b \psi \stackrel{def}{=} \sup_{ 1 \le p < b } \left[ \frac{|\xi|_p}{\psi(p)} \right],
 \  b = \const > 1,  \eqno(1.4)
$$
but in (1.4) $ \psi = \psi(p) $ is  continuous  function in the semi-open interval $ 1 \le p < b $
such that $ \lim_{p \uparrow b} \psi(p) = \infty; $  the case when $ \psi(b - 0) < \infty $ is trivial. \par

 \ We define formally in the case when $ b < \infty \ \psi(p) := + \infty \ $  for all the values $ p > b.$ \par
 \ An used further example:

 $$
 \psi^{(\beta,b)}(p) \stackrel{def}{=} (b-p)^{-\beta}, \ 1 \le p < b, \beta = \const > 0; \ G_{\beta,b}(p) := G_b\psi^{(\beta,b)}(p).
 $$

\vspace{3mm}

{\bf C.} Recall that sets $ A_1, A_2, \ A_i \in \cal{B} $ are named {\it disjoint,}  if $ A_1\cap A_2 = \emptyset. $ The sequence of functions
$  \{h_n \}, n =1,2,3  \ldots  $ is said to be {\it disjoint}, or more exactly {\it pairwise disjoint, if }

$$
\forall i,j: \ i \ne j \ \Rightarrow \ h_i \cdot h_j \stackrel{a.e.}{=} 0. \eqno(1.5)
$$

 \ If the sequence of (measurable) functions $  \{h_n \} $  is pairwise disjoint, then

$$
|\sum_n h_n|_p^p = \sum_n |h_n|_p^p, \hspace{6mm}  \sup_n |h_n(x)|^p = \sum_n |h_n(x)|^p, \ p = \const > 0. \eqno(1.6)
$$

\vspace{3mm}

{\bf D.} We denote as ordinary  for any measurable set $ A, \ A \in \cal{B} $ its indicator function by
$ I(A) = I_A(\omega).   $\par

\vspace{3mm}

{\bf E.} Let $ \xi = \xi(t), \ t \in T  $ be again separable random field (process) such that

$$
\exists b = \const \in(1, \infty], \ \forall p \in[1,b) \Rightarrow  \sup_{t \in T} |\xi(t)|_p < \infty.
$$
 \ Then the r.f. $ \xi(\cdot) $ generated  the so-called {\it natural } $ G\psi  \ - $ function by the formula

$$
\psi(p) = \psi^{(\xi)}(p) \stackrel{def}{=} \sup_{t \in T} |\xi(t)|_p, \ 1 \le p < b. \eqno(1.7)
$$

 \ Evidently,

$$
\forall t \in T \Rightarrow \xi(t) \in  G\psi^{(\xi)} \ {\bf and \ moreover} \    \sup_{t \in T} ||\xi(t)||G\psi^{(\xi)} = 1. \eqno(1.8)
$$

\vspace{3mm}

{\bf F.} Let $ \psi = \psi(p) $ be some function from the class $ G\psi_b, \ b = \const \in (1, \infty], $  such that
all the values $  \xi(t), \ t \in T $ belongs uniformly to the space $  G\psi. $ One can suppose without loss of generality

$$
\sup_{t \in T} ||\xi(t)||G\psi = 1. \eqno(1.9)
$$

 \ For instance, as the capacity of the function $ \psi(\cdot) $ may be picked the natural function  for the r.f. $  \xi: \
 \psi(p):=\psi^{(\xi)}(p), $ if of course there exists. \par

 \ Define by means of the function $  \psi(\cdot) $ the so - called {\it natural} (bounded) distance (more precisely, semi-distance)
$  d_{\psi}(t,s), \ t,s \in T $ on the set $  T: $

$$
d_{\psi}(t,s):= ||\xi(t) - \xi(s)||G\psi, \eqno(1.10)
$$
so that $   d_{\psi}(t,s) \le 2.   $\par
 \ Denote also by  $ D = D_{\psi} $ the diameter of the set $  T  $ relative the distance $ d_{\psi}:  $

$$
D =D(T,d_{\psi}) = \diam(T,d_{\psi}) \stackrel{def}{=} \sup_{t,s \in T} d_{\psi}(t,s), \eqno(1.11)
$$
 and by  $  H(T,d_{\psi},\epsilon)  $ the metric entropy of the set $  T  $  relative the distance $  d_{\psi} $ at the point
$ \epsilon, \ 0 < \epsilon < D.  $ \par
 \ The so-called entropy integral has by definition a form

$$
\Theta(T,d_{\psi}, \delta) \stackrel{def}{=} 9 \int_0^{\delta} \exp \left[ v_*( \ln 2 + H(T,d_{\psi},\epsilon)) \right] \ d \epsilon,
\ 0 < \delta \le D, \eqno(1.12)
$$
so that

$$
\Theta(T,d_{\psi}, D) \stackrel{def}{=} 9 \int_0^D \exp \left[ v_*( \ln 2 + H(T,d_{\psi},\epsilon)) \right] \ d \epsilon, \eqno(1.12a)
$$
where

$$
v(y) = \ln \psi(1/y), \  y \in (1/b,1); \hspace{5mm}  v_*(x) \stackrel{def}{=} \inf_{y \in (1/b,1)}(xy + v(y)). \eqno(1.13)
$$

 \ The transform $ v \to v_* $ is named as co-transform of Young-Fenchel, or Legendre, in contradiction to the classical
 Young-Fenchel transform

$$
f^*(x) \stackrel{def}{=} \sup_{y \in \dom f} (x y - f(y)). \eqno(1.14)
$$

 \ It is proved in particular in the monograph \cite{Ostrovsky3}, page 172, theorem 3.17.1, that if for the considered r.f.
$ \Theta(T,d_{\psi}, D) < \infty, $ then $ \xi(t) $ is $ d_{\psi}(\cdot, \cdot) $ continuous almost surely and herewith

$$
|| \ \sup_{t \in T} |\xi(t)| \ ||G\psi \le \Theta(T,d_{\psi}, D), \eqno(1.15)
$$
and following

$$
{\bf P} \left( \sup_{t \in T} |\xi(t)| > u \right) \le
\exp \left\{ - \nu^*_{\psi} \left[ \ \ln \left(u / \Theta(T,d_{\psi}, D) \right) \ \right] \right \},
\ u \ge \ 0. \eqno(1.16)
$$

 \ Moreover,

$$
|| \ \sup_{t,s: d_{\psi}(t,s) \le \delta} |\xi(t) - \xi(s)|  \ ||G\psi \le
$$

$$
 9 \int_0^{\delta} \exp \left[ v_*( \ln 2 + H(T,d_{\psi},\epsilon)) \right] \ d \epsilon = \Theta(T,d_{\psi}, \delta). \eqno(1.17)
$$

 \ The value

$$
\omega_{\xi, G\psi, d_{\psi}}(\delta) \stackrel{def}{=}   || \ \sup_{t,s: d_{\psi} \le \delta} |\xi(t) - \xi(s)|  \ ||G\psi
$$
is named as ordinary {\it probabilistic } module of continuity for the (uniform continuous) r.f. $ \xi(\cdot) $  relative the
distance function $  d_{\psi}. $  Obviously,

$$
 \lim_{\delta \downarrow 0} \omega_{\xi, G\psi, d_{\psi}}(\delta) = 0 \ \Rightarrow   {\bf P}_c(\xi) = 1.
$$

 \ Analogous estimates holds true is the r.f. $  \xi = \xi(t) $ satisfies the so-called majorizing measure condition, in particular, if
$ \sup_{t \in T} ||\xi(t)||G\psi < \infty,  $  see \cite{Talagrand1}, \ \cite{Talagrand2}-\cite{Talagrand5}.\par

\vspace{4mm}

\section{An interest example.}

\vspace{3mm}

 \ The following hypothesis has been formulated in the article \cite{Ostrovsky101}, 2008 year: \\

 \ "Let  $ \theta = \theta(t), t \in T  $  be arbitrary separable
random field, centered: $  {\bf E} \theta(t) = 0, $  bounded with probability one:
$ \sup_{t \in T} |\theta(t)| < \infty \ (\mod {\bf P }), $ moreover, may be continuous, if  the set $ T $
is   compact metric space relative some distance. \par
Assume in addition that for some Young (or Young-Orlicz) function $  \Phi(\cdot) $ and correspondent Orlicz norm $ ||\cdot||Or(\Phi) $

$$
\sup_{t \in T}||\theta(t)||Or(\Phi) < \infty. \eqno(2.1)
$$

 Recall that the Luxemburg norm $ ||\xi||Or(\Phi)  $ of a r.v. (measurable function) $  \xi $  is defined as follows:

 $$
  ||\xi||Or(\Phi) = \inf_{k, k > 0} \left\{ \int_{\Omega} \Phi(|\xi(\omega)|/k) \ {\bf P}(d \omega) \le 1 \right\}.
 $$

 \  The Young function $ \Phi(\cdot) $ is by definition arbitrary even convex continuous strictly increasing on the
non-negative right-hand  semi-axis  such that
$$
\Phi(0) = 0, \ \lim_{u \to \infty} \Phi(u) = \infty.
$$

 \ Let also $ \Psi(\cdot) $ be {\it arbitrary} another Young function such that $ \lim_{u \to \infty} \Psi(u) = \infty, \ \Psi <<  \Phi,  $
which denotes by definition

$$
\forall \lambda > 0  \ \Rightarrow \lim_{u \to \infty} \frac{\Psi(\lambda u)}{\Phi(u)} = 0,  \eqno(2.2)
$$
see \cite{Rao1}, p.16. \par
  \ Recall that $ \Psi << \Phi $ implies in particular that the (unit) ball in the space $ Or(\Psi) $ is precompact set
in the space $ Or(\Phi). $ \par

\vspace{4mm}

{\it Open question: there holds (or not)}"

$$
|| \sup_{t \in T} \ |\theta(t)| \ ||Or(\Psi) < \infty. \eqno(2.3)
$$

\vspace{4mm}

 \ A more slight question: there holds (or not)

$$
|| \sup_{t \in T} \ |\theta(t)| \ ||Or(\Phi) < \infty. \eqno(2.3a)
$$

 \ One can formulate the analogous question replacing the Orlicz spaces by Grand Lebesgue ones, see  \cite{Ostrovsky101}. \par

 \ The detail investigation of the theory of Orlicz's spaces including the case of unbounded source measure $  {\bf P} $ may be found
in the monographs \cite{Rao1}, \cite{Rao2}.\par

\vspace{4mm}

 \ This conclusions are true for the centered (separable) Gaussian fields,  \cite{Fernique1},
if the field $ \theta(\cdot) $
satisfies the so-called entropy or generic chaining condition \cite{Ostrovsky1},
\cite{Ostrovsky3},  \cite{Ostrovsky2}, \cite{Talagrand1}, \cite{Talagrand2}, \cite{Talagrand3},
\cite{Talagrand4}; in the case when  $ \theta(\cdot) $ belongs to the domain of attraction of
Law of Iterated Logarithm \cite{Ostrovsky5} etc.\par

 \ Notice that if the field $ \theta(t) $ is continuous $ ( \mod {\bf P}) $  and satisfies the condition
(2.1), then {\it there exists } an Young function $ \Psi(\cdot), \ \Psi(\cdot) << \Phi(\cdot) $
for which the inequality (2.3) there  holds, see \cite{Ostrovsky1}. \par

\vspace{3mm}

  \ The negative answers on these questions are obtained: in the article \cite{Ostrovsky100}, the case of Orlicz spaces;
in  \cite{Ostrovsky101}, more general and more strictly case of Grand Lebesgue Spaces. We recall briefly using further the
 correspondent example from the last report  \cite{Ostrovsky101}. \par

\vspace{4mm}

\ {\bf Example 2.1.} \par

\vspace{3mm}

{\bf 1.} We choose in the sequel in this pilcrow as  the capacity of compact metric space $  (T,d) $ the set of positive integer
numbers with infinite associated point which we denote by $  \infty: $

$$
T = \{ 1,2,3, \ldots, \infty   \}.  \eqno(2.4)
$$

 \ The distance $  d $ is defined as follows:

 $$
 d(i,j) = \left| \frac{1}{i} - \frac{1}{j} \right|, \ i,j < \infty; \ d(i,\infty) = d(\infty,i)= \frac{1}{i}, \ i < \infty;
 \eqno(2.5)
 $$
and obviously $ d(\infty,\infty) = 0.  $\par
 \ The pair $ (T,d) $ is compact (closed) metric space and the set $ T $ has an unique limit point $ t_{0} = \infty. $
For instance, $ \lim_{n \to \infty} d(n,\infty) = 0.   $ \par

\vspace{4mm}

{\bf 2.}  Let $ \Omega = (0,1)  $ with ordinary Lebesgue measure.
 \ Let  also  $  f = f(x), \ x \in \Omega = (0,1)  $ be non-zero non-negative integrable function belonging to the
space $  L_6.  $  Define a following $  \psi \ -  $ function:

$$
\nu(p) = |f|_p = \left[  \int_0^1 |f(x)|^p \ dx  \right]^{1/p}, \ 1 \le p \le 4.  \eqno(2.6)
$$
 \ On the other words, $  \nu(\cdot) $ is a natural function for the function $  f. $ Evidently,
$  \nu(\cdot) \in G\psi_4 = G\psi_{(0,4)}.  $\par

  \ Introduce also  the following numerical sequences

$$
c_n := n^{\beta},  \ \beta = \const > 0, \ n = 2,3, \ldots; \eqno(2.7)
$$

$$
\Delta_n := C(\beta) \cdot n^{ - 4 \beta - 1 }, \ C(\beta): \ \sum_{n=1}^{\infty} \Delta_n = 1; \
a_n = a(n):= \sum_{m=n}^{\infty} \Delta_n;  \eqno(2.8)
$$
and define also sequence of  functions and likewise the following  positive random process
$ \theta(t) = g_n, \ n = t, \ t,n \in T, \ \Omega = \{x \},  $

$$
g_n(x) = c(n) \ f \left( \frac{x-a(n)}{\Delta(n)} \right) \ I_{(a(n+1), a(n)) }(x), \ x \in \Omega, \
g_{\infty}(x) = 0;   \eqno(2.9)
$$

$$
g(x) = \sum_{n=1}^{\infty} g_n(x) =  \sum_{n=1}^{\infty} c_n  \ f \left( \frac{x-a(n)}{\Delta(n)} \right) \ I_{(a(n+1), a(n)) }(x).  \eqno(2.10)
$$

 \ Note that the sequence of r.v.  $ \{ g_n(x) \}  $ consists on the non-negative and disjoint functions, therefore

 $$
 \sup_n g_n(x) = \sum_n g_n(x) = g(x), \hspace{5mm} |\sup_n g_n|^p_p = \sum_n |g_n|^p_p. \eqno(2.11)
 $$

\vspace{3mm}

  \ Note also that the functions  $ g_n $ are disjoint and following $ \sup_n |g_n(x)| < \infty  $ almost surely. \par
We calculate using the relations (2.7)-(2.11):

$$
|g_n|_p^p = c^p(n) \ \Delta_n \ \nu^p(p) = C(\beta) \ n^{p \beta - 4 \beta - 1} \ \nu^p(p),
\ 1 \le p \le 4,  \eqno(2.12)
$$
therefore

$$
\sup_{p \in [1,4]}  \sup_n |g_n|_p^p \le C(\beta) \ \nu^4(4)  < \infty \eqno(2.13a)
$$
or equivalently

$$
\sup_n |g_n(\cdot)|_4 < \infty. \eqno(2.13b)
$$

 \ Moreover, $ g_n \to 0  $ almost everywhere.  Indeed, let $ \epsilon $ be arbitrary positive number. We get
applying the estimate (2.12) at the value $ p = 1 $  and Tchebychev-Markov inequality

$$
\sum_n {\bf P} ( |g_n| > \epsilon ) \le C(\beta) \sum_n \frac{n^{ - 3 \beta - 1  }}{\epsilon} < \infty.
$$
 Our conclusion follows immediately from the lemma of Borel-Cantelli. \par

\vspace{3mm}

 \  So, the random process $ \theta(t) = g_n, $ where $ \ n = t \ $ satisfies the condition (2.1)
relative the $ \Psi \ -   $ function $ \psi_{(4)}(p) := 1, \ 1 \le p \le 4  $ and is continuous almost everywhere relative the
source distance function $  d = d(t,s). $ \par

 \vspace{3mm}

\  Let us now find  the exact up to multiplicative constant expression for the natural function of the r.v.
 $  \sup_n |g_n(x)| $  as $  p \to 4-0.  $  We have:

 $$
  | \ \sup_n |g_n| \ |_p^p =  \sum_n |g_n|_p^p    = \sum_n c^p(n) \ \Delta_n \  \nu^p(p) =
 $$

$$
  = C(\beta) \ \ \nu^p(p)  \ \sum_n n^{p \beta - 4 \beta - 1} \ \sim \frac{C_1(\beta)}{4 - p};
  \eqno(2.14)
$$

$$
| \ \sup_n \ | g_n|  \ |_p  \sim C_2(\beta)(4 - p)^{-1/4}. \eqno(2.15)
$$

\vspace{3mm}

 \  Thus, we can choose for our  proposition  as capacity  the $ \Psi \ - $  function $ \psi(p) $  the
function  $  \psi_{(4)}(p), $ which is in turn equivalent to the following $   \Psi \ -  $ function

 $$
 \psi_{(4)}(p) := 1, \ 1 \le p < 4,
 $$
and correspondingly  to take

$$
\phi_0(p):=  (4 - p)^{-1/8} \stackrel{def}{=} \psi^{(1/8, 4)}(p), \ 1 \le p < 4. \eqno(2.16)
$$

 \ Obviously,

 $$
  \phi_0(\cdot) << \psi^{(0,4)}(\cdot) \eqno(2.17)
 $$
and

$$
|| \ \sup_n \ | g_n|  \ ||G\phi_0 = \infty, \eqno(2.18)
$$

\vspace{3mm}

 \ {\bf Remark 2.0.}  \ In order to obtain the centered needed process $ \theta(t)  $  with at the same properties,
 we consider the sequence  $  \tilde{g}_n(x) = \epsilon(n) \cdot g_n(x),  $  where $ \{\epsilon(n) \} $ is a Rademacher
 sequence independent on the $  \{g_n\} $  defined perhaps on some sufficiently rich probability space:

$$
{\bf P}(\epsilon(n) = 1) = {\bf P}(\epsilon(n) = -1) = 1/2; \eqno(2.19)
$$
then

$$
 | \ \tilde{g}_n(x) \ | = | \ g_n(x) \ |,  \ | \ \tilde{g}_n \ |_p = | \ g_n \ |_p \eqno(2.20)
$$
 and the sequence $ \{\tilde{g}_n \} $ is also pairwise disjoint (Rademacher's symmetrization).\par
 \ This completes the grounding of using for us properties of our (counter \ - ) example. \par

\vspace{4mm}

{\bf Remark 2.1.} The constructed process $ \theta(t) $ give us a new example of centered continuous random process with relatively
 light tails of finite-dimensional distribution, but for which the so-called entropy and generic chains series divergent.\par

\vspace{3mm}

{\bf Remark 2.2.} The properties of our example remains true if we use instead the space of continuous function $ C(T,d)  $
arbitrary  separable Banach space.\par

\vspace{4mm}

\section{Partition scheme. Main result-boundedness.}

\vspace{3mm}

 \ {\bf Definition 3.1.} The representation of the form

$$
T = \cup_{m=1}^{\infty} T_m, \eqno(3.1)
$$
write $  T \sim \{ T_m \}, $ on the (measurable) subsets $ \ T_m, \ $ not necessary to be disjoint,
is said to be a {\it partition, } or equally {\it covering} of the set $ T. $ \par

 \ We have for any partition

$$
{\bf P}_T(u) \le \sum_{m=1}^{\infty} {\bf P}_{T_m}(u), \ u > 0. \eqno(3.2)
$$

 \ In order to estimate each summands in (3.2), we need to introduce an "individual" natural parameters for the r.f. $ \xi(\cdot) $
on the arbitrary subset $  T_m. $ In detail, denote: \par

$$
\psi_m(p) := \sup_{t \in T_m} | \ \xi(t) \ |_p, \eqno(3.3)
$$
and assume

$$
\exists b_m \in (1, \infty], \ \forall p < b_m \ \Rightarrow \psi_m(p) < \infty. \eqno(3.4)
$$

 \ Define as before formally in the case $ b_m < \infty  \ \forall p > b_m \ \Rightarrow \psi(p) = + \infty. $ \par
 \ Further, put

$$
 v^{(m)}(y) := \ln \psi_m(1/y), \ y \in (0,1);    \ D_m := \diam(T_m, d_{\psi_m}), \eqno(3.5)
$$

$$
Z(m) :=  \Theta_m(T_m,d(\psi_m), D_m), \eqno(3.6)
$$

$$
Y( \{T_m\},u) \stackrel{def}{=} \sum_{m=1}^{\infty} \exp \left\{ - \nu^*_{\psi_m} \left[ \ \ln (u/ Z(m)) \ \right] \ \right\}. \eqno(3.7)
$$

\vspace{3mm}

{\bf Theorem 3.1.} Suppose that for some partition $   T \sim \{T_m\}  $

$$
\lim_{u \to \infty} Y( \{T_m\},u) = 0. \eqno(3.8)
$$
 Then the r.f. $ \xi(t) $ is bounded with probability one and moreover

$$
{\bf P} ( \sup_{t \in T} |\xi(t)| > u) \le Y( \{T_m\},u). \eqno(3.9)
$$

\vspace{3mm}

 \ {\bf Proof.} We have using (1.16)

$$
{\bf P} \left( \sup_{t \in T_m} |\xi(t)| > u \right) \le \exp \left\{ - \nu^*_{\psi_m} \left[ \ \ln (u/ Z(m)) \ \right] \ \right\},
\ u \ge \ 0. \eqno(3.10)
$$

 \ The proposition (3.9) follows immediately from the estimate (3.2). \par

\vspace{3mm}

 \ {\bf Corollary 3.1.}  We conclude under conditions of theorem (3.1)

$$
{\bf P} ( \sup_{t \in T} |\xi(t)| > u) \le \inf_{ \{T_m\}  } Y( \{T_m\},u)  \stackrel{def}{=}R(u), \ u \ge 0. \eqno(3.11)
$$
where $ "\inf" $ in (3.11) is calculated over all the partitions $ \ T \sim \{ T_m \} $ of the set $  T. $ \par

\vspace{3mm}

 \ {\bf Remark 3.1.} The estimates (3.10) and (3.11) are exponential non-improvable still  in the previous entropy approach,
 i.e. without partition, see \cite{Ostrovsky3}, chapters 3,4.\par

\vspace{3mm}

  {\bf  An example. } It is easy to see that the our example  in the second section satisfies all the conditions of theorem 3.1. The
 partition for the considered therein r.f. is trivial. \par

\vspace{4mm}

\section{Partition scheme. Main result-continuity.}

\vspace{3mm}

 \ The condition (3.8) guarantee us only the {\it boundedness} of almost all the paths of the r.f. $ \xi(t). $ We will discuss
 in this section the sufficient conditions relative the appropriate  distance
 for the {\it continuity} of $ \xi(\cdot) $  based on the offered in this article partition scheme. \par
 \ We will follow the article of V.A.Dmitrovsky \cite{Dmitrovsky1}.\par

 \ Some new notations.  Introduce a new $ \psi \ - $ function $  \tau = \tau(p) $ as follows:

$$
\tau(p) := \left[ p \ \int_0^{\infty} u^{p-1} \ R(u) \ du \right]^{1/p}, \eqno(4.1)
$$
so that

$$
\left[{\bf E} \left\{ \sup_{t \in T} |\xi(t)| \right\}^p \right]^{1/p} \le \tau(p), \eqno(4.2)
$$
and suppose its finiteness at last for one value $ p \ > 1; $ define $  b = \sup\{p, \ p > 1, \tau(p) < \infty; $
then $ \sup_t ||\xi(t)|| G\tau \le 1.  $\par

\ The correspondent bounded natural distance $ \rho = \rho(t,s) $ may be defined as follows:

$$
 \rho(t,s) := ||\xi(t) - \xi(s)||G\tau. \eqno(4.3)
$$

 \ We find as before using (1.17)

$$
|| \ \sup_{t,s: \rho(t,s) \le \delta} |\xi(t) - \xi(s)|  \ ||G\tau  \le  \Theta(T,\rho, \delta), \eqno(4.4)
$$
or equally
$$
\omega_{\xi, G\tau, \rho}(\delta) \le \Theta(T,\rho, \delta), \ \delta \in (0, \diam(T, \rho)). \eqno(4.5)
$$

\ Thus, we proved in fact the following statement. \par

\vspace{4mm}

\ {\bf  Proposition 4.1. } We conclude under formulated above in this section notations and conditions:
 the r.f. $ \xi = \xi(t) $ is $ \rho(\cdot, \cdot) $ is uniform continuous with probability one. \par

\vspace{4mm}

 \ {\bf  Example 4.1 - 2.1.} Let us return to the example 2.1 considered in the second section. We deduce after some calculations
taking into account proposition 4.1

$$
 {\bf P}(  \ \sup_n \ | g_n| > u  ) \le \frac{C_1(\beta) \ \ln u}{u^4}, \ u > e. \eqno(4.6)
$$
 \ The continuity a.e. of the r.p. $ g_n $ is proved in the second section; we intend to prove further the probabilistic
continuity. It is sufficient to consider the unique limit point $  n = \infty. $ \par

 \ Define the next Young function

$$
\Phi(u) = e^{2} \ u^2, \ |u| \le e, \hspace{5mm} \Phi(u) = \frac{u^4}{\ln |u|}, \ |u| > e.
$$

 \ The distance function $ \rho(n,m) $ is here equivalent to the Orlicz's natural distance

$$
r(n,m) := || \ g_n - g_m \ ||L(\Phi)
$$
and is equivalent in turn to the source distance $ d(m,n), \ m,n = 1,2,\ldots, \infty.  $ \par
 \ Our statement:

$$
\lim_{n \to \infty} || \ \sup_{k \ge n} |g_k| \ ||L(\Phi) = 0,
$$
 by virtue of proposition 4.1. \par


\vspace{4mm}

\begin{thebibliography}{99}

\vspace{3mm}

\bibitem{Dudley1}
{\sc Dudley R.M.} {\it The sizes of compact of Hilbert space and continuity of Gaussian
processes.} J. Functional Analysis, (1967), B. 1 pp. 290-330.

\bibitem{Dmitrovsky1}
 {\sc Dmitrovsky V.A.} (1981). {\it  On the distributions of maximum and local properties of realizations
 of the pre-gaussian fields.}  Theory of Probab. and Math. Stat.; (in Russian). Kiev, KSU, {\bf 25}, 154-164.

\bibitem{Fernique1}
{\sc Fernique X.} (1975). {\it Regularite des trajectoires des function aleatiores
gaussiennes.} Ecole de Probablite de Saint-Flour, IV, Lecture Notes in Mathematics;
480, {\bf 1,} 96, Springer Verlag, Berlin.

\bibitem{Fernique2}
{\sc Fernique X,} {\it Caracterisation de processus de trajectoires majores ou
continues.} Seminaire de Probabilites XII. Lecture Notes in Math.,  649,
(1978), 691–706, Springer, Berlin.

\bibitem{Fernique3}
{\sc Fernique X.} {\it Regularite de fonctions aleatoires non gaussiennes.} Ecolee
 de Ete de Probabilit'es de Saint-Flour XI,  Lecture Notes in Math.,
 976, (1983), 1–74, Springer, Berlin.

\bibitem{Fiorenza1}
   {\sc Capone C., Fiorenza A., Krbec M.} {\it On the Extrapolation Blowups in the
   $ L_p $ Scale. } Collectanea Mathematica, {\bf 48}, 2, (1998), 71-88.

 \bibitem{Fiorenza2}
 {\sc Fiorenza A.} {\it Duality and reflexivity in grand Lebesgue spaces.}
       Collectanea Mathematica (electronic version), {\bf 51}, 2, (2000), 131-148.

\bibitem{Fiorenza3}
 {\sc Fiorenza A., and Karadzhov G.E.} {\it Grand and small Lebesgue spaces and
       their analogs.} Consiglio Nationale Delle Ricerche, Instituto per le
      Applicazioni del Calcoto Mauro Picine, Sezione di Napoli, Rapporto tecnico n. 272/03, (2005).

\bibitem{Frolov1}
{\sc Frolov A.S., Tchentzov N.N. } {\it On the calculation by the Monte-Carlo
method definite integrals depending on the parameters. } Journal of Computational
Mathematics and Mathematical Physics, (1962), V. 2, Issue 4, p. 714-718 (in Russian).

\bibitem{Grigorjeva1}
{\sc Grigorjeva M.L., Ostrovsky E.I.} {\it Calculation of Integrals on discontinuous
Functions by means of depending trials method.} Journal of Computational
Mathematics and Mathematical Physics, (1996), V. 36, Issue 12, p. 28-39 (in Russian).

\bibitem{Iwaniec1}
   {\sc Iwaniec T., and Sbordone C.} {\it On the integrability of the Jacobian under
      minimal hypotheses.} Arch. Rat.Mech. Anal., 119, (1992), 129–143.

\bibitem{Iwaniec2}
 {\sc Iwaniec T., P. Koskela P., and Onninen J.} {\it Mapping of finite distortion:
   Monotonicity and Continuity. } Invent. Math. 144 (2001), 507-531.

\bibitem{Kolmogorov1}
{\sc Kolmogorov A.N.} {\it About analitical methods in the probability  theory.} Surveys of Soviet Math., (1938),
B.17, issue 5, 5-41.

\bibitem{Kozachenko1}
 {\sc Kozachenko Yu. V., Ostrovsky E.I.} (1985). {\it The Banach Spaces of
      random Variables of subgaussian type.}  Theory of Probab. and Math.
      Stat.; (in Russian). Kiev, KSU, {\bf 32}, 43-57.

\bibitem{Kurbanmuradov1}
{\sc Kurbanmuradov O., Sabelfeld K.} (2007). {\it Exponential bounds for
   the probability deviation of sums of random fields.} Preprint.
   Weierstra$\beta$ - Institut f\"ur Angewandte Analysis und Stochastik
   (WIAS), ISSN 0946-8633, p. 1-16.

\bibitem{Talagrand1}
{\sc Ledoux M., Talagrand M.} (1991) {\it Probability in Banach Spaces.}
      Springer, Berlin, MR 1102015.

\bibitem{Liflyand1}
{\sc Liflyand E., Ostrovsky E., Sirota L.} {\it Structural Properties of Bilateral Grand Lebesgue Spaces.}
Turk. Journal of Math., 34, (2010), 207-219.

 \bibitem{Ostrovsky1}
 {\sc Ostrovsky E.} {\it Support of Borelian Measures in separable Banach Spaces.} \\
 arXiv:0808.03248v1 [math.FA] 24 Aug 2008.

\bibitem{Ostrovsky2}
{\sc Ostrovsky E., Rogover E. } {\it Exact exponential Bounds for the random field Maximum Distribution
via the majorizing Measures (Generic Chaining).} \\
arXiv:0802.0349v1 [math.PR] 4 Feb 2008.

\bibitem{Ostrovsky3}
{\sc  Ostrovsky E.I.} (1999). {\it Exponential estimations for random Fields and its
applications}, (in Russian). Moscow-Obninsk, OINPE.

\bibitem{Ostrovsky5}
{\sc Ostrovsky E.I.} (1994.) {\it Exponential Bounds  in the Law of Iterated
Logarithm in Banach Spaces.} Math. Notes, {\bf 56}, 5, p. 98-107.

\bibitem{Ostrovsky6}
{\sc Ostrovsky E., Sirota L.} {\it Fourier Transforms in exponential rearrangement invaiant Spaces.}
arXiv:040639v1 [math.FA], 20 Jun 2004.

\bibitem{Ostrovsky7}
{\sc Ostrovsky E., Rogover E.} {\it Maximal Inequalities in bilateral Grand Lebesgur Spaces over unbounded Measure.} \\
arXiv:0808v1 [math.FA], 24 Aug 2008.

\bibitem{Ostrovsky8}
{\sc  Ostrovsky E. and Sirota L.} {\it Moment Banach spaces: theory and applications.}
HIAT Journal of Science and Engineering, {\bf C}, Volume 4, Issues 1-2, pp. 233-262, (2007).

 \bibitem{Ostrovsky9}
 {\sc  Ostrovsky E. and Sirota L.} {\it Moment and Tail Inequalities for polynomial Martimgales.
 The case of heavy tails.}\\
 arXiv: 1112.2768v1 [math.PR] 13 Dez 2011.

\bibitem{Ostrovsky100}
{\sc  Ostrovsky E. and Sirota L.} {\it A counterexample to a hypothesis of light tail of maximum distribution for continuous
random processes with light  finite-dimensional tails. } \\
arXiv:1208.6281v1 [math.PR] 30 Aug 2012

\bibitem{Ostrovsky101}
{\sc  Ostrovsky E. and Sirota L.} {\it  Can the tail for maximum of continuous random field be significantly more heavy
 than maximum of tails?  }\\
arXiv:1508.05646v1 [math.PR] 23 Aug 2015


\bibitem{Pizier1}
{\sc Pizier G.} {\it Condition  d' \ entropic  assupant la continuite de certains processus et
applications a l’analyse harmonique.} Seminaire  d \ analyse fonctionalle. (1980),
Exp.13, p. 23-34.

\bibitem{Rao1}
{\sc Rao M.M., Ren Z.D.} {\it Theory of Orlicz Spaces.}  Marcel Dekker Inc., 1991. New
York, Basel, Hong Kong.

\bibitem{Rao2}
{\sc Rao M.M., Ren Z.D.} {\it Applications of Orlicz Spaces.} Marcel Dekker Inc., 2002.
New York, Basel, Hong Kong.

\bibitem{Slutzky1}
{\sc Slutzky E.E.} {\it Some propositions about the theory of random functions.}  Proceedings of Middle-Asia
University,   Tashkent,  Ser. Math., (5), {\bf 31,} (1949), 3-15. (in Russian).

 \bibitem{Talagrand2}
 {\sc Talagrand M.} (1996). {\it Majorizing measure: The generic chaining.}
      {\it Ann. Probab.} {\bf 24} 1049-1103. MR1825156

\bibitem{Talagrand3}
 {\sc Talagrand M.} (2001). {\it Majorizing Measures without Measures.}
    {\it Ann. Probab.} 29, 411-417. MR1825156

\bibitem{Talagrand4}
 {\sc Talagrand M.} (2005).  {\it The Generic Chaining. Upper and
     Lower Bounds of Stochastic Processes.} Springer, Berlin. MR2133757.

\bibitem{Talagrand5}
  {\sc Talagrand M.}(1990). {\it Sample boundedness of stochastic processes
      under increment conditions.}  Ann. Probab., {\bf 18}, 1-49.

\vspace{4mm}

\end{thebibliography}
\end{document}